\documentclass[12pt, oneside, reqno]{amsart}

\usepackage{amsmath}
\usepackage{amsfonts}
\usepackage{tikz}
\usetikzlibrary{calc}
\usepackage{graphicx}
\usepackage{bm}
\usepackage{geometry}
\usepackage{lscape}
\usepackage{caption}
\usepackage{subcaption}
\usepackage{url}
\usepackage{hyperref}
\usepackage{setspace}
\usepackage[numbers]{natbib}
\usepackage{mathtools}
\usepackage{cancel}
\usepackage{soul}
\usepackage{dutchcal}
\usepackage{helvet}
\allowdisplaybreaks 
\usepackage{algorithm}
\usepackage{algpseudocode}

\title[Efficient implementation of the hybridized Raviart-Thomas method]{Efficient implementation of the hybridized Raviart-Thomas mixed
method by converting \\ flux subspaces into stabilizations}

\author[S. Anantharamu]{Sreevatsa Anantharamu}
\address[S. Anantharamu]{Senior Applications Engineer at X-ScaleSolutions,
Minneapolis, Minnesota, USA}
\email{anant035@umn.edu} 

\author[B. Cockburn]{Bernardo Cockburn}
\address[B. Cockburn]{School of Mathematics, University of Minnesota, Minneapolis, MN 55455, USA}
\email{bcockbur@umn.edu}

\thanks{Bernardo Cockburn's research was supported in part by the Advanced Computational Center for Entry Systems Simulation (ACCESS) through NASA grant 80NSSC21K1117.}

\subjclass[2010]{Primary 65N30, 65M60, 35L65}

\keywords{mixed methods, hybridization, static condensation, mixed methods, hybridizable discontinuous Galerkin methods}

\date\today

\begin{document}

\newcommand{\mbf}[1]{\mathbf{#1}} 
\newcommand{\what}[1]{\hat{#1}}
\newcommand{\pder}[2]{\frac{\partial #1}{\partial #2}} 
\newcommand{\pderev}[3]{\frac{\partial #1}{\partial #2}\rvert_{#3}} 
\newcommand{\ppder}[3]{\frac{\partial^2 #1}{\partial #2 \partial #3}} 
\newcommand{\psqder}[2]{\frac{\partial^2 #1}{\partial {#2}^2}}
\newcommand{\psqderev}[3]{\frac{\partial^2 #1}{\partial {#2}^2}\rvert_{#3}}
\newcommand{\pa}{\partial}
\newcommand{\larr}{\leftarrow}
\newcommand{\rarr}{\rightarrow}
\newcommand{\mapping}[3]{#1 : #2 \rightarrow #3}

\newcommand{\ip}[2]{\left( #1 , #2 \right)} 
\newcommand{\ipthr}[3]{\left( #1 #2 , #3 \right)} 
\newcommand{\ipOm}[2]{\left( #1 , #2 \right)_{\Omega}} 
\newcommand{\ipK}[2]{( #1 , #2 )_{K}} 
\newcommand{\ipdom}[3]{\left( #1 , #2 \right)_{#3}} 
\newcommand{\crb}[1]{\{#1\}} 
\newcommand{\bob}[1]{\left[#1\right]} 
\newcommand{\subh}[1]{{#1}_\mathrm{h}}
\newcommand{\subhK}[1]{{#1}_{\mathrm{h},K}}
\newcommand{\subhKg}[1]{{#1}_{\mathrm{h},K,\bm{g}}}
\newcommand{\subhKlhat}[1]{{#1}_{\mathrm{h},K,\what{\lambda}_{\mathrm{h}}}}
\newcommand{\subhKmu}[1]{{#1}_{\mathrm{h},K,\mu}}
\newcommand{\subhKpsii}[1]{{#1}_{\mathrm{h},K,\psi_i}}
\newcommand{\subhKpsij}[1]{{#1}_{\mathrm{h},K,\psi_j}}
\newcommand{\subhKp}[1]{{#1}_{\mathrm{h},K^+}}
\newcommand{\subhKm}[1]{{#1}_{\mathrm{h},K^-}}

\newcommand{\ipv}[2]{\left(\mathbf{#1},\mathbf{#2}} 
\newcommand{\ipvOm}[2]{\left(\mathbf{#1} , #2 \right)_{\Omega}} 
\newcommand{\ipvK}[2]{\left( #1 , #2 \right)_{K}} 
\newcommand{\ipvdom}[3]{\left( #1 , #2 \right)_{#3}} 
\newcommand{\dprfunc}[2]{#1 \cdot #2} 
\newcommand{\ipbdom}[3]{\langle #1, #2\rangle_{#3}} 

\newcommand{\wtilde}[1]{\widetilde{#1}}
\newcommand{\wtil}[1]{\widetilde{#1}}
\newcommand{\bvphi}{\bm{\varphi}}
\newcommand{\itbPhi}{\bm{\mathit{\Phi}}}
\newcommand{\bxi}{\bm{\xi}}
\newcommand{\Omegahat}{\what{\Omega}}
\newcommand{\bnabla}{\bm{\nabla}}
\newcommand{\bu}{\bm{u}}
\newcommand{\butil}{\wtilde{\bu}}
\newcommand{\buhat}{\what{\bm{u}}}
\newcommand{\buK}{\mathbf{u}_K}
\newcommand{\lam}{\lambda}
\newcommand{\lamK}{\lambda_K}
\newcommand{\bU}{\bm{U}}
\newcommand{\bv}{\bm{v}}
\newcommand{\bV}{\bm{V}}
\newcommand{\bL}{\bm{L}}
\newcommand{\bn}{\bm{n}}
\newcommand{\mT}{\mathrm{T}}
\newcommand{\bvg}{\bm{v}_{\mathrm{g}}}
\newcommand{\bVg}{\bm{V}_{\mathrm{g}}}
\newcommand{\mrm}[1]{\mathrm{#1}}
\newcommand{\trm}[1]{\textrm{#1}}
\newcommand{\FinvT}{\bm{F}^{-\mrm{T}}}
\newcommand{\Finvv}{\bm{F}^{-1}}
\newcommand{\Fninvv}{\bm{F}^{n^{-1}}}
\newcommand{\FninvT}{\bm{F}^{n^{-\mrm{T}}}}
\newcommand{\Fnpinvv}{\bm{F}^{{n+1}^{-1}}}
\newcommand{\FnpinvT}{\bm{F}^{{n+1}^{-\mrm{T}}}}
\newcommand{\Omegahatijk}{\Omegahat_{i,j,k}}
\newcommand{\ijk}[1]{{#1}_{i,j,k}}
\newcommand{\rst}[1]{{#1}_{r,s,t}}
\newcommand{\mcal}[1]{\mathcal{#1}}
\newcommand{\ijkl}[1]{{#1}_{i,j,k,l}}
\newcommand{\ovbar}[1]{\overline{#1}}
\newcommand{\bmcal}[1]{\bm{\mcal{#1}}}
\newcommand{\bcdot}{\boldsymbol{\cdot}}
\newcommand{\ijkv}[2]{{#1}_{i,j,k,#2}}
\newcommand{\bvpsi}{\bm{\psi}}
\newcommand{\invT}[1]{{#1}^{-\mrm{T}}}
\newcommand{\inv}[1]{{#1}^{-1}}
\newcommand{\bX}{\bm{X}}
\newcommand{\bx}{\bm{x}}
\newcommand{\gradn}[1]{\left[G_{\mcal{N}}#1\right]}
\newcommand{\ijkx}[2]{{#1}_{i,j,k,#2}}
\newcommand{\GbXf}[2]{\ijkx{\left[\bm{G}_{\bX}#1\right]}{#2}}
\newcommand{\pmatthree}[3]{\begin{pmatrix}#1\\#2\\#3\end{pmatrix}}
\newcommand{\Gammahat}{\what{\Gamma}}
\newcommand{\gsbt}[2]{#1_{\!_{#2}}} 
\newcommand{\bnablaX}{\gsbt{\bnabla}{\bX}}
\newcommand{\DX}{\gsbt{\bm{D}}{\bX}}
\newcommand{\mrmf}{\mrm{f}}
\newcommand{\mrms}{\mrm{s}}
\newcommand{\mrmfs}{\mrmf\mrms}
\newcommand{\bus}{\bu_{\mrms}}
\newcommand{\buf}{\bu_{\mrmf}}
\newcommand{\bUs}{\bU_{\mrms}}
\newcommand{\bUf}{\bU_{\mrmf}}
\newcommand{\bvphif}{\bvphi_{\trm{f}}}
\newcommand{\bvphis}{\bvphi_{\trm{s}}}
\newcommand{\Omegaf}{\Omega_{\trm{f}}}
\newcommand{\Omegas}{\Omega_{\trm{s}}}
\newcommand{\mrmn}{\mrm{n}}
\newcommand{\ticv}{\trm{icv}}
\newcommand{\trmh}{\trm{h}}
\newcommand{\qtil}{\wtilde{q}}
\newcommand{\Qtil}{\wtilde{Q}}
\newcommand{\xhat}{\what{x}}
\newcommand{\yhat}{\what{y}}
\newcommand{\zhat}{\what{z}}
\newcommand{\bftil}{\wtilde{\bm{f}}}
\newcommand{\bLtil}{\wtilde{\bm{L}}}
\newcommand{\pmpbar}{\left(p-\ovbar{p}\right)}
\newcommand{\bphihat}{\bm{\what{\phi}}}
\newcommand{\mrmuhat}{\what{\mrm{u}}}

\newcommand{\red}[1]{{\color{red}#1}}
\newcommand{\blah}{SOMETHING??? }

\newcommand{\Pkd}{\left[P_k\right]^d}
\newcommand{\Pk}{P_k}
\newcommand{\Ptilk}{\wtilde{P}_k}
\newcommand{\Ptil}[1]{\wtilde{P}_{#1}}
\newcommand{\Ptilkd}{\left[\wtilde{P}_k\right]^d}
\newcommand{\Ptild}[1]{\left[\wtilde{P}_{#1}\right]^d}
\newcommand{\RR}{\mathbb{R}}
\newcommand{\PkK}{P^k(K)}
\newcommand{\PkOm}{P^k(\Omega)}
\newcommand{\Pkdom}[1]{P^k(#1)}
\newcommand{\Ltwo}{L^2}
\newcommand{\Hone}{H^1}
\newcommand{\Om}{\Omega}
\newcommand{\OmK}{\Omega_K}
\newcommand{\paOm}{\partial \Omega}
\newcommand{\paOmK}{\partial \Omega_K}
\newcommand{\Vh}{V_h}
\newcommand{\union} {\cup}
\newcommand{\belongs}{\in}
\newcommand{\grad}{\nabla}
\newcommand{\eps}{\epsilon}
\newcommand{\xhatvec}{\vec{\what{x}}}

\newcommand{\myint}[3]{\int_{#1} #2\,\mrmd\Omega}
\newcommand{\ipiK}[2]{\int_K #1 #2\,d\Omega}
\newcommand{\ipi}[3]{\int_{#3} #1 #2\,d\Omega}
\newcommand{\ipivK}[2]{\int_K {#1} \cdot {#2}\,d\Omega}
\newcommand{\ipiv}[3]{\int_{#3} {#1} \cdot {#2}\,d\Omega}
\newcommand{\ipibdom}[3]{\int_{#3}#1#2\,d\Gamma}

\newcommand{\jump}[1]{\llbracket #1 \rrbracket}

\newcommand{\expnumber}[2]{{#1}\mathrm{e}{#2}}
\newcommand{\enu}[2]{\expnumber{#1}{#2}}

\newcommand{\dunderline}[1]{\underline{\underline{#1}}}
\newcommand{\bmu}{\bm{\mu}}
\newcommand{\qhat}{\what{q}}
\newcommand{\bpsi}{\bm{\psi}}
\newcommand{\mrmd}{\mrm{d}}
\newcommand{\onebytwo}{\frac{1}{2}}
\newcommand{\threebytwo}{\frac{3}{2}}
\newcommand{\rdots}{\red{\dots}}
\newcommand{\nmhalf}{n-\onebytwo}
\newcommand{\bg}{\bm{g}}
\newcommand{\phihat}{\hat{\phi}}
\newcommand{\avg}[1]{\{\!\{#1\}\!\}}
\newcommand{\uhat}{\what{u}}

 \begin{abstract}
     We show how to reduce the computational time of the practical implementation of the Raviart-Thomas mixed method for second-order elliptic problems. The implementation takes advantage of a recent result which states that certain local subspaces of the vector unknown can be eliminated from the equations by transforming them into stabilization functions; see the paper published online in JJIAM on August 10, 2023. We describe in detail the new implementation (in MATLAB and a laptop with Intel(R) Core (TM) i7-8700 processor which has six cores and hyperthreading) and present numerical results showing 10 to 20\% reduction in the computational time for the Raviart-Thomas method of index $k$, with $k$ ranging from 1 to 20, applied to a model problem.
     \end{abstract}
     \maketitle
 
\section{Introduction}
Let us begin by noting that the title of the paper contains both words: `mixed methods' and `stabilization'. At first glance, this
might appear to be a technical error since mixed methods do not need
stabilization since stability is directly ensured by the choice of their function spaces. 
However, this is \textit{not} an error. In a recent paper,
\cite{bcockburn2023hybridizable}, one of the authors showed how a portion
of the flux space of any hybridized mixed method can be recast as a stabilization of a equivalent hybridizable Discontinuous Galerkin (HDG) methods. By
{\em equivalent}, we mean that the original hybridized mixed method and the
new HDG methods result in exactly the same solution. 
In this paper, we show that the implementation of the {\em equivalent} HDG method is faster than that of original
hybridized mixed method. We carry out this for {\em two}
equivalent HDG methods.

The paper is organized as follows. In Section 2, we describe how to obtain an HDG method from a hybridized mixed method; we follow \cite{bcockburn2023hybridizable}.
In Section 3, we describe in full detail the implementation of the hybridized Raviart-Thomas method. We then do the same for {\em two} equivalent HDG methods.
In Section 4, we display our numerical results. We end with some concluding remarks.  We use the standard notation used for HDG methods, see, for example, \cite{CockburnEncy18,bcockburn2023hybridizable}.


\section{Background}

\subsection{Subspace-to-stabilization.} For the sake of completeness, we begin by summarizing the subspace-to-stabilization
result \cite[Section 5]{bcockburn2023hybridizable}. Consider
the Poisson problem:
\begin{equation} \label{eqn:poisson problem}
  \begin{split}
    {c}\,\bm{q} =& -\bnabla u, \text{ in }\Omega,\text{ and}\\
    \bnabla \cdot \bm{q} =& f,\text{ in }\Omega,
  \end{split}
\end{equation}
with the boundary condition that $u=u_D$ on the boundary $\partial
\Omega$. Here $f$ is the source term and $u_D$ is the Dirichlet boundary
data. The hybridized mixed method formulation for this problem is as
follows: For each element $K$ of the mesh, find $\bm{q}_h$ and $u_h$
belonging to the function spaces $\bm{V}$ and $W$, respectively, such
that:
\begin{equation} \label{eqn:local problem hrt}
  \begin{split}
    \ipK{c\,\bm{q}_h}{\bm{v}} - \ipK{u_h}{\bnabla\cdot\bv} =& -\ipbdom{\hat{u}_h}{\bv\cdot\bn}{\partial K},\\
    \ipK{\bnabla\cdot \bm{q}_h}{w} =& \ipK{f}{w},
  \end{split}
\end{equation}
for all test functions $\bm{v}$ and $w$ belonging to the function spaces
$\bm{V}$ and $W$, respectively. The above equations are referred to as
the `local problem'. Here, $\what{u}_h$ is an approximation to $u$ on
the faces of the triangulation and is a data to the above problem. The
additional equations for this face variable are:
\begin{equation}\label{eqn:global problem hrt}
  \begin{split}
    \ipbdom{\bm{q}_h^+\cdot\bn^+ + \bm{q}_h^-\cdot\bn^-}{\mu}{F_I} =& 0 , \text{ for all $\mu\in M_h(F_I)$}, \\
    \ipbdom{\what{u}_h}{\mu}{F_D} =& \ipbdom{u_D}{\mu}{F_D} , \text{ for all $\mu\in M_h(F_D)$}, \\
  \end{split}
\end{equation}
where $F_I$ is each interior face of the mesh, and $F_D$ is each
Dirichlet boundary face of the mesh. The above equations are referred to
as the `global problem'.

Split $\bm{V}$ into $\bm{V}_s \oplus \bm{V}_a$. Here, $\bm{V}_s$ is the
subspace of $\bm{V}$ that will be converted into stabilization and
$\bm{V}_a$ is the subspace of $\bm{V}$ that will be used to define the
local problem of the equivalent HDG method. Their exact definition is
deferred until later. This splitting converts the local problem
(equation \eqref{eqn:local problem hrt}) into:
\begin{equation*}
  \begin{split}
    \ipK{c \left(\bm{q}_a + \bm{q}_s\right)}{\bm{v}_a} - \ipK{u_h}{\bnabla\cdot\bv_a} =& -\ipbdom{\hat{u}_h}{\bv_a\cdot\bn}{\partial K},\\
    \ipK{c \left(\bm{q}_a + \bm{q}_s\right)}{\bm{v}_s} - \ipK{u_h}{\bnabla\cdot\bv_s} =& -\ipbdom{\hat{u}_h}{\bv_s\cdot\bn}{\partial K},\\
    \ipK{\bnabla\cdot \left( \bm{q}_a + \bm{q}_s \right)}{w} =& \ipK{f}{w},
  \end{split}
\end{equation*}
for all test functions $\bm{v}_a$, $\bm{v}_s$, and $w$ in the function
spaces $\bm{V}_a$, $\bm{V}_s$, and $W$, respectively. Requiring the
functions in $\bm{V}_a$ and $\bm{V}_s$ to be orthogonal to each other in
the $(c\cdot,\cdot)_K$ inner-product, we obtain:
\begin{equation}\label{eqn:split local prob}
  \begin{split}
    \ipK{c\, \bm{q}_a}{\bm{v}_a} - \ipK{u_h}{\bnabla\cdot\bv_a} =& -\ipbdom{\hat{u}_h}{\bv_a\cdot\bn}{\partial K},\\
    \ipK{c\, \bm{q}_s}{\bm{v}_s} - \ipK{u_h}{\bnabla\cdot\bv_s} =& -\ipbdom{\hat{u}_h}{\bv_s\cdot\bn}{\partial K},\\
    \ipK{\bnabla\cdot \left( \bm{q}_a + \bm{q}_s \right)}{w} =& \ipK{f}{w}.
  \end{split}
\end{equation}
Integrating the second equation by parts and requiring $\bm{V}_s$ to be
any subspace of $\bm{V}$ that is $L^2$-orthogonal to $\bnabla W$ yields
\begin{equation}\label{eqn:qs}
  \ipK{c\, \bm{q}_s}{\bm{v}_s} = \ipbdom{u_h-\hat{u}_h}{\bv_s\cdot\bn}{\partial K}.
\end{equation}
Note that in the above equation, $\bm{q}_s$ (which is the portion of $q$
in the subspace $\bm{V}_s$) depends solely on the jump $u_h-\what{u}_h$
on the faces of the element. Observe that the appearance of the term
$u_h-\what{u}_h$ has some similarity to the flux stabilization
$\tau(u_h-\what{u}_h)$ that is used in a HDG method (note that
$(\tau(\cdot)$ here is a linear-mapping that satisfies certain necessary
requirements). This similarity is exploited to define a stabilization
function $\tau(u_h-\what{u}_h)$ for the HDG method that is equivalent to
the above hybridized mixed method below.

Based on the form of equation \eqref{eqn:qs}, let us define
$\bm{L}_{\bm{V}_s}^c$ to be the lifting operator that maps a function
$\mu$ in the space $W$ to the function $\bm{L}_{\bm{V}_s}^c(\mu)$ in the
space $\bm{V}_s$ as:
\begin{equation} \label{eqn: eq stab}
  \ipK{c\, \bm{L}_{\bm{V}_s}^c(\mu)}{\bm{v}_s} = \ipbdom{\mu}{\bv_s\cdot\bn}{\partial K},
\end{equation}
for all test functions $\bm{v}_s$ in the space $\bm{V}_s$. Then,
$\bm{q}_s=\bm{L}_{\bm{V}_s}^c(u_h-\what{u}_h)$. Moreover, the hybridized
mixed method in equation \eqref{eqn:split local prob} can be manipulated to
the following HDG method for the portion $\bm{q}_a$ of the flux
approximation and the full scalar approximation $u_h$: Find $\bm{q}_a$
and $u_h$ in the function spaces $\bm{V}_a$ and $W$, respectively, such
that:
\begin{equation}\label{eqn:eqv hdg method local}
  \begin{split}
    \ipK{c\, \bm{q}_a}{\bm{v}_a} - \ipK{u_h}{\bnabla\cdot\bv_a} &= -\ipbdom{\hat{u}_h}{\bv_a\cdot\bn}{\partial K},\\
    -\ipK{\bm{q}_a}{\bnabla w} + \ipbdom{\bm{q}_h\cdot\bm{n}}{w}{\partial K}&= \ipK{f}{w}, \text{and}\\
    \bm{q}_h\cdot\bm{n} &= \bm{q}_a\cdot\bm{n} + \tau(u_h-\what{u}_h),
  \end{split}
\end{equation}
where the stabilization function $\tau(\cdot)$ is defined using the
lifting operator as: $\tau(u_h-\what{u}_h)$ =
$\bn\cdot\bm{L}_{\bm{V}_s}^c(u_h-\what{u}_h)$. The global problem of the
HDG method is same as the one in equation \eqref{eqn:global problem hrt}.

The flux approximation $\bm{q}_h$ of the hybridized mixed method in
equation \eqref{eqn:local problem hrt} is then obtained from the above HDG
method as $\bm{q}_h=\bm{q}_a+\bm{L}_{\bm{V}_s}^c(u_h-\what{u}_h)$. Thus,
the effect of the space $\bm{V}_s$ is fully encapsulated in the defined
stabilization function $\tau$ via the lifting operator
$\bm{L}_{\bm{V}_s}^c$. This is the crux of the `spaces-to-stabilization'
idea.

In summary, the conditions on the spaces $\bm{V}_a$ and $\bm{V}_s$ are that:
\begin{equation}\label{eqn:spaces to stabilization conds}
    \begin{aligned}
        \text{$\bm{V}_s$ should be any subspace of $\bm{V}$ that is $L^2$-orthogonal to $\nabla W$,\phantom{ool}}\\
        \text{$\bm{V}_a$ should be the remaining portion of $\bm{V}$ that is orthogonal to $\bm{V}_s$}\\
        \text{in the $(c\,\cdot,\cdot)_K$ inner-product. 
        \phantom{oooooooooooooooooooooooooolooo}}
    \end{aligned}
\end{equation}
The former and latter conditions above were used to obtain the equations
\eqref{eqn:qs} and \eqref{eqn:split local prob}, respectively.

\subsection{The equivalent HDG method after static condensation}. Similar to
other HDG methods, the degrees of freedom corresponding to $\bm{q}_a$
and $u_h$ can be statically condensed to yield a globally-coupled
problem just for the degrees of freedom corresponding to the face
variable $\what{u}_h$ as follows. The local problem for the equivalent
HDG method given in equation \eqref{eqn:eqv hdg method local} can be shown to
be equal to the following local problem: 

\

\noindent Find $ \bm{q}_a$ and $u_h$ in
the function spaces $\bm{V}_a$ and $W$, respectively, such that
\begin{equation}\label{eqn:eq hdg method local1}
  \begin{split}
    \ipK{c\, \bm{q}_a}{\bm{v}_a} - \ipK{u_h}{\bnabla\cdot\bv_a} =& -\ipbdom{\hat{u}_h}{\bv_a\cdot\bn}{\partial K},\\
    \ipK{\bnabla \cdot \bm{q}_a}{w} + \ipK{c \bm{L}_{\bm{V}_s}^c \left(u_h\right)}{\bm{L}_{\bm{V}_s}^c \left(w\right)} &= \ipK{f}{w} + \ipK{c \bm{L}_{\bm{V}_s}^c \left( \what{u}_h \right)}{\bm{L}_{\bm{V}_s}^c \left(w\right)},
  \end{split}
\end{equation}
for all test functions $\bm{v}_a$ and $w$ in the function spaces
$\bm{V}_a$ and $W$, respectively. The influence of $\what{u}_h$ and $f$
on $(\bm{q}_a,u_h)$ can be separated as $(\bm{q}_a,u_h)$ =
$(\bm{Q}_{\what{u}_h},U_{\what{u}_h})$ + $(\bm{Q}_f,U_f)$. Here,
$(\bm{Q}_{\what{u}_h},U_{\what{u}_h})$ is the solution to the following
problem with $\mu=\what{u}_h$: 

\

\noindent Find $\bm{Q}_{\mu}$ and $U_{\mu}$ in the
function spaces $\bm{V}_a$ and $W$, respectively, such that:
\begin{equation}\label{eqn:eq hdg method local uhat}
  \begin{split}
    \ipK{c\,\bm{Q}_{\mu}}{\bm{v}_a} - \ipK{U_{\mu}}{\bnabla\cdot\bv_a} =& -\ipbdom{\mu}{\bv_a\cdot\bn}{\partial K},\\
    \ipK{\bnabla \cdot \bm{Q}_{\mu}}{w} + \ipK{c \bm{L}_{\bm{V}_s}^c \left( U_{\mu} \right)}{\bm{L}_{\bm{V}_s}^c \left(w\right)} &= \ipK{c \bm{L}_{\bm{V}_s}^c \left( \mu \right)}{\bm{L}_{\bm{V}_s}^c \left(w\right)},
  \end{split}
\end{equation}
for all test functions $\bm{v}_a$ and $w$ in the function spaces
$\bm{V}_a$ and $W$, respectively. $(\bm{Q}_f,U_f)$ is the solution to
the problem:

\

\noindent Find  $\bm{Q}_f$ and $U_f$ in the function spaces
$\bm{V}_a$ and $W$, respectively, such that:
\begin{equation}\label{eqn:eq hdg method local f}
  \begin{split}
    \ipK{c\,\bm{Q}_{f}}{\bm{v}_a} - \ipK{U_{f}}{\bnabla\cdot\bv_a} =& 0\\
    \ipK{\bnabla \cdot \bm{Q}_{f}}{w} + \ipK{c \bm{L}_{\bm{V}_s}^c \left( U_{f} \right)}{\bm{L}_{\bm{V}_s}^c \left(w\right)} &= \ipK{f}{w},
  \end{split}
\end{equation}
for all test functions $\bm{v}_a$ and $w$ in the function spaces
$\bm{V}_a$ and $W$, respectively. 

Using the above decomposition, we can
show that the equation for $\what{u}_h$ given in equation
\eqref{eqn:global problem hrt} is equal to the following  problem: 

\

\noindent Find
$\what{u}_h$ belonging to the function space $M_h$ such that:
\begin{equation}\label{eqn: eq hdg method global}
{\footnotesize
\begin{aligned}
        \sum_K \ipK{c\, \bm{Q}_{\uhat_h}}{\bm{Q}_{\mu}} 
        + \sum_K \ipK{c\, \bm{L}_{\bm{V}_s}^c (U_{\uhat_h} - \uhat_h)}{ \bm{L}_{\bm{V}_s}^c (U_{\mu}-\mu)} =& \sum_K \ipK{f}{U_{\mu}}, \quad\forall \mu\in M_h(F_I),\\
    \ipbdom{\what{u}_h}{\mu}{F_D} =& \ipbdom{u_D}{\mu}{F_D}, \phantom{oooo}\forall \mu\in M_h(F_D),
    \end{aligned}
}
\end{equation}
for each interior face $F_I$ and Dirichlet boundary face $F_D$ of the
triangulation.

\section{Implementations of the Hybridized Raviart-Thomas mixed
method}

We use this subspace-to-stabilization idea to come up with two new
implementations of the hybridized Raviart-Thomas (RT) mixed method.
Each implementation stems from a different choice of the subspace
$\bm{V}_s$.  We note that for the hybridized RT method, the local spaces $\bm{V}$ and $W$
are: 
\begin{equation*}
   \bm{V} = [P_k(K)]^d \oplus \bm{x}\wtilde{P}_k(K)\text{ and } W = P_k(K). 
\end{equation*}
Here, ${P}_k(K)$ denotes the space polynomials of degree $k$ defined on $K$, whereas  $\wtilde{P}_k(K)$ denotes the space of homogeneous polynomials of degree $k$ defined on $K$. Table \ref{table:choices of Vs and Va} shows the three choices of $\bm{V}_s$ and
the name of the implementation that stems from each of these choices.

\begin{table}[h!]
    \centering
    \begin{tabular}{c|c|c|c} 
     \hline
      Implementation  & $\bm{V}_s$  & Notes & $\bm{V}_a$  \\ 
     \hline\hline
      Usual-HRT   & $\bm{V}_s^{(0)} (=\{\bm{0}\})$ & Usual & $[P_k(K)]^d +
      \bm{x} P_k(K) $ \\
      Stab-1-HRT  & $\bm{V}_s^{(1)}$           & New   & $[P_k(K)]^d$ \\ 
      Stab-2-HRT  & $\bm{V}_s^{(2)}$           & New   & $[P_{k-1}(K)]^d$ \\
     \hline
    \end{tabular}
    \caption{The different implementations of the hybridized RT method. Note that $\bm{V}^{(1)}_s$ is the largest subspace of $\bm{V}$ that is $\bm L^2$-orthogonal to $[P_k(K)]^d$, and that $\bm{V}^{(2)}_s$
    is of the form $\bm{V}_s^{(1)}\oplus\bm{V}^{(3)}$  where
    $\bm V^{(3)}$ is the subspace of  polynomials in $[P_k(K)]^d$ that are
      $\bm L^2$-orthogonal to $[P_{k-1}(K)]^d$.}
    \label{table:choices of Vs and Va}
\end{table}

The
implementation Usual-HRT is the usual implementation of the hybridized RT
method, see \cite{ArnoldBrezzi85}. This stems from the choice $\bm{V}_s=\bm{V}_s^{(0)}=\{\bm{0}\}$ (the
empty set). The implementations Stab-1-HRT and Stab-2-HRT are the two new
implementations proposed in this paper. The new implementation Stab-1-HRT stems
from the choice $\bm{V}_s=\bm{V}_s^{(1)}$, where $\bm{V}_s^{(1)}$ is the largest 
subspace of $\bm{V}$ that is $L^2-$orthogonal to $[P_k(K)]^d$. The other new
implementation $\bm{V}_s^{(2)}$ stems from choosing $\bm{V}_s=\bm{V}_s^{(2)}$,
where $\bm{V}_s^{(2)}$ is the space $\bm{V}_s^{(2)}$ plus the vector-valued
polynomials in $[P_k(K)]^d$ that are $L^2$-orthogonal to $[P_{k-1}(K)]^d$.


For each of these implementations, the local space $\bm{V}_a$ is the
$(c\cdot,\cdot)_K$-orthogonal complement of $\bm{V}_s$ within $\bm{V}$ .
In this paper, we will consider the case $c=Id$. Hence, the
$(c\cdot,\cdot)_K$ inner-product becomes the standard $L^2$
inner-product. Therefore, the space $\bm{V}_a$ for the Usual-HRT,
Stab-1-HRT, and Stab-2-HRT implementation is $[P_k(K)]^d+\bm{x}P_k(K)$,
$[P_k(K)]^d$, and $[P_{k-1}(K)]^d$, respectively. The details of each of these 
implementations are given next.

\clearpage
\subsection{Usual-HRT}

The details of the Usual-HRT implementation are given below.

\subsubsection{Basis} \label{subsubsec: basis usual-hrt}

In the Usual-HRT implementation, the space $\bm{V}_a$ equals 
$[P_k(K)]^d \oplus \bm{x}\wtilde{P}_k(K)$.
We use the following basis functions for this space in each element $K$: 
\begin{equation*}
  \bvphi_1^{(K)},\dots,\bvphi_{n}^{(K)}.
\end{equation*}
Here, $n=dm + m'$, where $m=C_{d}^{k+d-1}$ and $m'=C_{d-1}^{k+d-1}$. The
first $dm$ basis functions $\bvphi_1^{(K)},\dots,\bvphi^{(K)}_{dm}$ correspond to
the $[P_k(K)]^d$ portion of the space. The remaining $m'$ functions
$\bvphi_{dm+1}^{(K)}, \dots, \bvphi_{n}^{(K)}$ correspond to the remaining portion
of the space. 

These basis functions are orthonormal to each other. They satisfy the orthonormality condition:
\begin{equation} \label{eqn:bvphi orth}
  \int_{K}\bvphi_i^{(K)}\cdot\bvphi_j^{(K)}\,\mrmd\Omega = |K|\delta_{ij},
\end{equation}
where $\delta_{ij}$ is the Kroenecker delta and $|K|$ is the measure 
(area in 2D and volume in 3D) of element $K$.
The first $dm$ basis functions $\bvphi_1,\dots,\bvphi_{dm}$
are constructed using the (normalized) Dubiner polynomials
\citep{dubiner1991spectral} as:
\begin{equation} \label{eqn:bvphi first dm}
  \bvphi_{d(i'-1)+j'}^{(K)} = q_{i'}^{(K)}\bm{e}_{j'}, \text{ for $i'=1,m$ and
  $j'=1,\dots,d$.}
\end{equation}
Here, $q_1,\dots,q_m$ are the orthonormal Dubiner polynomials in the
element $K$. These are obtained by mapping the orthonormal polynomials
in the reference simplex $\what{K}$ to the element $K$ as: 
\begin{equation*}
  q_i^{(K)}(\bm{x}^{(K)}(\what{\bm{x}})) = \what{q}_i(\what{\bm{x}}),
\end{equation*}
where $\bm{x}^{(K)}(\what{\bm{x}})$ is the affine mapping from the
reference element $\what{K}$ to the element $K$. 
These polynomials satisfy the following orthonormality relation:
\begin{equation} \label{eqn:q orth}
  \int_K q_i^{(K)} q_j^{(K)} \, \mrmd\Omega = |K| \delta_{ij}.
\end{equation}
The remaining $m'$
basis functions $\bvphi_{dm + 1}^{(K)},\dots,\bvphi_n^{(K)}$ are
constructed by multiplying the degree $k$ Dubiner polynomial basis
functions with $\bm{x}$ and orthonormalizing them with the rest of the
basis functions using the modified Gram-Schmidt kernel. This is described
in algorithm \ref{algo: addn RT basis}.

\begin{algorithm}
  \caption{Generating $\bvphi_{dm+1}^{(K)},\dots,\bvphi_n^{(K)}$}
  \begin{algorithmic}
    \For{$i=dm+1,\dots,n$}
      \State $\bvphi_i^{(K)} \gets \bm{x} q^{(K)}_{i-(d-1)m-m'}$
      \State Orthonormalize $\bvphi_i^{(K)}$ against the previous $(i-1)$ basis
      functions using a modified Gram-Schmidt kernel
    \EndFor
  \end{algorithmic}
  \label{algo: addn RT basis}
\end{algorithm}

For the space $W$ (which equals $P_k(K)$), we use the same orthonormal Dubiner 
polynomial basis functions:
\begin{equation*}
  q_1^{(K)},\dots,q_m^{(K)}.
\end{equation*}
For the space $M_h$, we also use orthonormal Dubiner polynomial basis functions but defined on
the faces of the mesh. They are:
\begin{equation*}
  \psi_1^{(F)},\dots,\psi_{m'}^{(F)},
\end{equation*}
where $F$ is a typical face of the mesh.




\subsubsection{Local problem} \label{subsubsec:loc prob usual hrt}

Using the above basis functions converts the local problem that depends on $\mu$ (equation
\eqref{eqn:eq hdg method local uhat}) to the following matrix problem:
\begin{equation}\label{eqn:lin alg local prob mu}
  \begin{bmatrix}
    \left[I_{n\times n}\right] & -\left[D_{m \times n}^{(K)}\right]^T\\
    \left[D^{(K)}_{m\times n}\right] & \left[0_{m\times m}\right]
  \end{bmatrix}
  \begin{bmatrix}
    \left[Q^{\mu,(K)}_{n\times ((d+1)m')}\right] \\
    \left[U^{\mu,(K)}_{m\times ((d+1)m')}\right]
  \end{bmatrix}=
  \begin{bmatrix}
    -\left[b^{Q^{\mu},(K)}_{n\times ((d+1)m')}\right]\\
    \left[0_{m\times ((d+1)m')}\right]
  \end{bmatrix}.
\end{equation}
Here, $\left[I_{n\times n}\right]$ is the identity matrix, $\left[0_{a\times b}\right]$ is 
$a\times b$ matrix of zeros,
$\left[D^{(K)}_{m\times n}\right]$ is the divergence matrix,
$\left[Q^{\mu,(K)}_{n\times (d+1)m'}\right]$ is the 
degree of freedom matrix corresponding to the element-wise mapping $\bm{Q}_{\mu}$, 
$\left[U^{\mu,(K)}_{m\times ((d+1)m')}\right]$ is the 
degree of freedom matrix corresponding to the element-wise mapping $U_{\mu}$, and 
$-\left[b^{Q^{\mu},(K)}_{n\times ((d+1)m')}\right]$ is the right-hand side matrix 
corresponding to the degrees of freedom 
$\left[Q^{\mu,(K)}_{n\times (d+1)m'}\right]$. The size of all matrices is given 
in the subscript. The entries of the divergence and right-hand matrices are:
\begin{equation}
    \begin{split}
        \left[D^{(K)}_{m\times n}\right]_{i,j}=\int_K q_i \bnabla \cdot \bvphi^{(K)}_j \, \mrmd\Omega, \text{ and }\left[b^{Q^{\mu},(K)}_{n\times ((d+1)m')}\right]_{i,(j-1)m' + r} = \int_{F_j}\psi_{r}^{(F_j)}\bvphi^{(K)}_i \cdot \bn \, \mrmd\Gamma,
    \end{split}
\end{equation}
respectively, where $F_j$ is the $j^{th}$ face of element $K$.


For efficient solution of the above local problem, we perform two optimizations. 
The first optimization pertains to computing the divergence matrix 
$\left[D^{(K)}_{m\times n}\right]$. This matrix has the following form:
\begin{equation} \label{eqn:div mat usual-hrt}
  \left[D^{(K)}_{m\times n}\right]=
  \begin{bmatrix}
    \begin{matrix}
      \left[D^{(K)^{(1)}}_{(m-m')\times dm}\right] \\
      \left[\phantom{o.o}0_{m'\times dm}\phantom{oo}\right]
    \end{matrix} & \left[D^{(K)^{(2)\phantom{^|}}}_{m\times m'{\phantom{_\big|}}}\right]
  \end{bmatrix}.
\end{equation}
Here, $\left[D^{(K)^{(1)}}_{(m-m')\times dm}\right]$ is the portion of the
divergence-matrix from the first $dm$ basis functions,
$\bvphi_1^{(K)},\dots,\bvphi_{dm}^{(K)}$, which correspond to the $[P_k(K)]^d$ portion of the
space. The remaining portion of the divergence matrix,
$\left[D^{(K)^{(2)}}_{m\times m'}\right]$, is from the last $m'$ basis
functions, $\bvphi_{dm+1}^{(K)},\dots,\bvphi_{n}^{(K)}$. This form is exploited for its efficient 
computation as follows. The portion
$\left[D^{(K)^{(1)}}_{(m-m')\times dm}\right]$ is first computed in along the 
coordinates of the reference element $\what{K}$. We
will denote this reference divergence matrix as $\left[\what{D}_{(m-m')\times dm}\right]$ and 
is given by:
\begin{equation*}
    \left[\what{D}_{(m-m')\times dm}\right]_{i,j}=\int_{\what{K}} \what{q}_i \what{\bnabla} \cdot \what{\bvphi}_j \, \mrmd\what{\Omega}.
\end{equation*}
Then, to compute
$\left[D^{(K)^{(1)}}_{(m-m')\times dm}\right]$ in each element $K$, we simply
combine the columns of $\left[\what{D}_{(m-m')\times dm}\right]$ using the entries of the 
Jacobian of $\bm{x}^{(K)}(\what{\bm{x}})$. Hence, we do not differentiate the first $dm$ 
basis functions separately for each element of the mesh but just once
in the reference element. Elementwise differentiation is performed only for the 
last $m'$ basis functions needed to compute the remaining of the divergence matrix,
portion, $\left[D^{(K)^{(2)}}_{m\times m'}\right]$.

The second optimization pertains to the solution of the local matrix problem 
in equation \eqref{eqn:lin alg local prob mu}. Since the mass matrix is identity, the 
degrees of freedom matrix $\left[Q^{\mu,(K)}_{n\times (d+1)m'}\right]$ can be easily
eliminated to obtain the following matrix problem just for the degrees of freedom 
$\left[U^{\mu,(K)}_{m\times ((d+1)m')}\right]$:
\begin{equation}\label{eqn:lap local prob mu}
  \left[L^{(K)}_{m\times m}\right]\left[U^{\mu,(K)}_{m\times ((d+1)m')}\right] = \left[D^{(K)}_{m\times n}\right] \left[b^{Q^{\mu},(K)}_{n\times ((d+1)m')}\right].
\end{equation}
Here, $\left[L^{(K)}_{m\times m}\right]$ is the Laplacian matrix and is equal to:
\begin{equation*}
    \left[L^{(K)}_{m\times m}\right]=\left[D^{(K)}_{m\times n}\right]
\left[D^{(K)}_{m\times n}\right]^T.
\end{equation*}
To solve this problem, the matrix $\left[L^{(K)}_{m\times m}\right]$ is first factored using the
(dense) Cholesky factorization method and the computed factor is used to obtain 
$\left[U^{\mu,(K)}_{m\times ((d+1)m')}\right]$.

Similarly, the above basis functions convert the local problem that depends on $f$ (equation \eqref{eqn:eq hdg method local f}) to the following matrix problem:
\begin{equation*}
   \begin{bmatrix}
    \left[I_{n\times n}\right] & \left[-D_{m \times n}^{(K)}\right]^T\\
    \left[D^{(K)}_{m\times n}\right] & \left[0_{m\times m}\right]
  \end{bmatrix}
  \begin{bmatrix}
    \left[Q^{f,(K)}_{n\times 1}\right] \\
    \left[U^{f,(K)}_{m\times 1}\right]
  \end{bmatrix}=
  \begin{bmatrix}
    \left[0_{n\times 1}\right]\\
    |K|\left[P^{(K)}f_{m\times 1}\right]
  \end{bmatrix}. 
\end{equation*}
Here, $\left[Q^{f,(K)}_{n\times 1}\right]$ and $\left[U^{f,(K)}_{m\times 1}\right]$
are the degrees of freedom vector corresponding to the mapping $\bm{Q}_f$ and $U_f$,
respectively.Here, $\left[P^{(K)}f_{m\times 1}\right]$ is the degree of freedom 
vector obtained by $L^2$-projection of $f$ onto $W$ using the above basis functions. 
Similar to the local problem that depended on $\mu$, the above matrix problem is also solved 
by eliminating the degrees of freedom corresponding to $\left[Q^{f,(K)}_{n\times 1}\right]$ 
and then using the above computed Cholesky factor to obtain $ \left[U^{f,(K)}_{m\times 1}\right]$.



\subsubsection{Global problem} \label{subsubsec: global problem usual-hrt}

Using the above computed local problem solutions, the element matrix 
$\left[A^{(K)}_{((d+1)m') \times ((d+1)m')}\right]$ 
and element vector $\left[b^{(K)}_{((d+1)m') \times 1}\right]$ are computed in 
each element $K$ as:
\begin{equation}\label{eqn:elem mat and vec}
\begin{split}
   \left[A^{(K)}_{((d+1)m') \times ((d+1)m')}\right] &= \left[Q^{\mu,(K)}_{n\times (d+1)m'}\right]^T \left[Q^{\mu,(K)}_{n\times (d+1)m'}\right] |K|,\text{ and } \\
   \left[b^{(K)}_{((d+1)m') \times 1}\right] &= \left[U^{\mu,(K)}_{m\times (d+1)m'}\right]^T \left[P^{(K)}f_{m\times 1}\right] |K|,
\end{split}
\end{equation}
respectively. These element 
matrices and vectors are assembled together and the degrees of freedom of $\what{u}_h$ 
corresponding to the Dirichlet boundary faces are statically condensed to obtain the 
global matrix problem:
\begin{equation}\label{eqn: global mat problem}
  \left[A_{m'n_F \times m'n_{F}}\right] \left[\what{\mrm{u}}_{m'n_{F} \times 1}\right]=
  \left[ b_{m'n_{F} \times 1} \right].
\end{equation}
Here, $n_F$ is the total number of faces of the mesh minus the Dirichlet boundary faces.
$\left[A_{m'n_F \times m'n_{F}}\right]$ is the global left-hand side matrix.
$\left[ b_{m'n_{F} \times 1} \right]$ is the global right-hand side vector.
$\left[\what{\mrm{u}}_{m'n_{F} \times 1}\right]$ is the degree of freedom vector
corresponding to $\what{u}_h$. The above global matrix problem is solved using the sparse
Cholesky factorization method.



\subsection{Stab-1-HRT (new)}

The details of the new Stab-1-HRT implementation are given below. For general HDG methods, see \cite{CockburnGopalakrishnanLazarov09,CockburnDurham16}.

\subsubsection{Basis}

In this implementation, for the space
$\bm{V}_a$ in each element $K$, we use the following orthonormal basis functions:
\begin{equation*}
  \bvphi_1^{(K)},\dots,\bvphi_{dm}^{(K)}.
\end{equation*}
The above functions $\bvphi_i^{(K)}$ are the same ones
that were defined in section \ref{subsubsec: basis usual-hrt}. 
For the space $\bm{V}_s$, we use:
\begin{equation*}
  \bvphi_{dm+1}^{(K)},\dots,\bvphi_n^{(K)}.
\end{equation*}
These are the remaining $m'$ basis functions 
that were defined in section \ref{subsubsec: basis usual-hrt}.
For the remaining spaces $W$ and $M_h$, we use 
the same basis functions as that used in section \ref{subsubsec: basis usual-hrt}.

\subsubsection{Local problem} \label{subsubsec: local problem stab-1-hrt}

In this implementation, we have to compute the stabilization mapping $\bm{L}_{\bm{V}_s}^c$ 
(given in equation \eqref{eqn: eq stab}). Using the above basis functions for $\bm{V}_s$, the
matrix form of this mapping becomes:
\begin{equation*}
  \left[L_{m'\times (d+1)m'}^{s,(K)}\right]=\frac{1}{|K|}\left[b^{s,(K)}_{m'\times (d+1)m'}\right].
\end{equation*}
Here, $\left[b^{s,(K)}_{m'\times (d+1)m'}\right]$ is the right-hand side of the stabilization 
mapping given by:
\begin{equation*}
    \left[b^{s,(K)}_{m'\times ((d+1)m')}\right]_{i,(j-1)m' + r} = \int_{F_j}\psi_{r}^{(F_j)}\bvphi^{(K)}_{dm+i} \cdot \bn \, \mrmd\Gamma.
\end{equation*}


Using the above basis functions for $\bm{V}_a$ and the above matrix form of the 
stabilization mapping, the local problem that depends on $\mu$ 
(equation \eqref{eqn:eq hdg method local uhat}) becomes the following local matrix problem:
\begin{equation}\label{eqn:lin alg local prob mu stab-1}
  \begin{bmatrix}
    \left[I_{n\times n}\right] & \left[-D_{m \times dm}^{(K)}\right]^T\\
    \left[D^{(K)}_{m\times dm}\right] & \left[M^{s,(K)}_{m \times m}\right]
  \end{bmatrix}
  \begin{bmatrix}
    \left[Q^{\mu,(K)}_{dm\times (d+1)m'}\right] \\
    \left[U^{\mu,(K)}_{m\times (d+1)m'}\right]
  \end{bmatrix}=
  \begin{bmatrix}
    -\left[b^{Q^{\mu},(K)}_{dm\times ((d+1)m')}\right]\\
    \left[b^{U^{\mu},(K)}_{m\times ((d+1)m')}\right]
  \end{bmatrix}
\end{equation}
Here, $\left[M^{s,(K)}_{m\times m}\right]$ is the mass matrix that arises from
the equivalent stabilization. It relates to the stabilization mapping matrix as:
\begin{equation*}
    \left[M^{s,(K)}_{m\times m}\right] = |K|\left[\left[L_{m'\times (d+1)m'}^{s,(K)}\right]B^{(K)}_{(d+1)m'\times m}\right]^T \left[\left[L_{m'\times (d+1)m'}^{s,(K)}\right]B^{(K)}_{(d+1)m'\times m}\right],
\end{equation*}
where $B^{(K)}_{(d+1)m'\times m}$ is the linear mapping from the degrees of freedom of 
$u$ to that of $\mu$ in element $K$. $D^{(K)}_{m\times n}$ is the same divergence matrix that was
defined in section \ref{subsubsec:loc prob usual hrt}. The right-hand side matrix 
$\left[b^{U^{\mu},(K)}_{m\times ((d+1)m')}\right]$ is given by:
\begin{equation*}
    \left[b^{U^{\mu},(K)}_{m\times ((d+1)m')}\right] =  |K|\left[\left[L_{m'\times (d+1)m'}^{s,(K)}\right]B^{(K)}_{(d+1)m'\times m}\right]^T \left[L_{m'\times (d+1)m'}^{s,(K)}\right].
\end{equation*}


Observe that the divergence matrix $D^{(K)}_{m\times dm}$ in 
equation \eqref{eqn:lin alg local prob mu stab-1} has the following form:
\begin{equation*}
  \left[D^{(K)}_{m\times dm}\right]=
  \begin{bmatrix}
    \begin{matrix}
      \left[D^{(K)^{(1)}}_{(m-m')\times dm}\right] \\
      \left[\phantom{o.o}0_{m'\times dm}\phantom{oo}\right]
    \end{matrix}
  \end{bmatrix},
\end{equation*}
where the portion $\left[D^{(K)^{(1)}}_{(m-m')\times dm}\right]$ is the same as that
defined in equation \eqref{eqn:div mat usual-hrt}. Observe that the portion
$\left[D^{(K)^{(2)}}_{n\times m'}\right]$ that was present in equation
\eqref{eqn:div mat usual-hrt} is not present above. Hence, we only compute the
portion $\left[D^{(K)^{(1)}}_{(m-m')\times dm}\right]$ using the optimized
implementation described in section \ref{subsubsec:loc prob usual hrt}. We 
never compute $\left[D^{(K)^{(2)}}_{n\times m'}\right]$ that required differentiation of the
last $m'$ basis functions separately in each element of the mesh. This yields in
considerable reduction in the time taken to compute the individual local
problems and is demonstrated by the numerical experiments in section \ref{sec: numerical results}.


Similar to the Usual-HRT implementation, the degree of freedom matrix
$\left[Q^{\mu,(K)}_{dm\times (d+1)m'}\right]$ is eliminated to obtain the following 
equation for $\left[U^{\mu,(K)}_{m\times (d+1)m'}\right]$:
\begin{equation}\label{eqn:lap local prob mu stab-hrt-1}
   \left[L^{(K),1}_{m\times m}\right] \left[U^{\mu,(K)}_{m\times (d+1)m'}\right] = \left[D^{(K)}_{m\times n}\right] \left[b^{Q^{\mu},(K)}_{dm\times ((d+1)m')}\right] + \left[b^{U^{\mu},(K)}_{m\times ((d+1)m')}\right].
\end{equation}
Here, the Laplacian matrix $\left[L^{(K),1}_{m\times m}\right]$ is given by:
\begin{equation*}
    \left[L^{(K),1}_{m\times m}\right] = \left[D^{(K)}_{m\times dm}\right]
\left[D^{(K)}_{m\times dm}\right]^T + \left[M_{m\times m} ^{s,(K)}\right]. 
\end{equation*}
The Laplacian is first factored using the
Cholesky factorization method and the computed factors are used to compute 
$\left[U^{\mu,(K)}_{m\times (d+1)m'}\right]$. 

Similarly, the chosen basis functions convert the local problem that depends on $f$ to
the following matrix problem:
\begin{equation*}
   \begin{bmatrix}
    \left[I_{n\times n}\right] & \left[-D_{m \times dm}^{(K)}\right]^T\\
    \left[D^{(K)}_{m\times dm}\right] & \left[M_{m\times m}^{s,(K)}\right]
  \end{bmatrix}
  \begin{bmatrix}
    \left[Q^{f,(K)}_{dm\times (d+1)m'}\right] \\
    \left[U^{f,(K)}_{m\times (d+1)m'}\right]
  \end{bmatrix}=
  \begin{bmatrix}
    \left[0_{dm\times 1}\right] \\
   |K|\left[P^{(K)}f_{m\times 1}\right]
  \end{bmatrix}. 
\end{equation*}
The above matrix problem is solved similar to the above local problem that depended only on
$\mu$. The same Cholesky factors are used for its solution.


\subsubsection{Global problem}

In this implementation, the element matrix and vector are computed as:
\begin{equation}\label{eqn:elem mat and vec stab-1-hrt}
\begin{split}
   \left[A^{(K)}_{((d+1)m') \times ((d+1)m')}\right] &= \left[Q^{\mu,(K)}_{dm\times (d+1)m'}\right]^T \left[Q^{\mu,(K)}_{dm\times (d+1)m'}\right] |K|,\text{ and } \\
   \left[b^{(K)}_{((d+1)m') \times 1}\right] &= \left[U^{\mu,(K)}_{m\times (d+1)m'}\right]^T \left[P^{(K)}f_{m\times 1}\right] |K|.
\end{split}
\end{equation}
Then, similar to the Usual-HRT implementation, the above element matrices and element vectors are assembled to form the global problem: 
\begin{equation*}
  \left[A_{m'n_F \times m'n_{F}}\right] \left[\what{\mrm{u}}_{m'n_{F} \times 1}\right]=
  \left[ b_{m'n_{F} \times 1} \right].
\end{equation*}
This global problem is solved using the sparse Cholesky factorization.

\subsection{Stab-2-HRT (new)}

The details of the new Stab-2-HRT implementation are given below. 

\subsubsection{Basis}

In this implementation, for the space $\bm{V}_a$ in each element $K$,
we use the following orthonormal basis functions:
\begin{equation*}
  \bvphi_1^{(K)},\dots,\bvphi_{d(m-m'))}^{(K)},
\end{equation*}
where the above functions $\bvphi_i^{(K)}$ are the same functions
defined in section \ref{subsubsec: basis usual-hrt}. 
For the space $\bm{V}_s$, we use the remaining $(d+1)m'$ functions as the basis:
\begin{equation*}
  \bvphi_{d(m-m')+1}^{(K)},\dots,\bvphi_n.
\end{equation*}
For the spaces $W$ and $M_h$, we use the same basis functions that were used in section \ref{subsubsec: basis usual-hrt}.

\subsubsection{Local problem}

The local problem solution for this implementation is very similar to that of the
Stab-1-HRT implementation. The
matrix form of the stabilization mapping $\bm{L}_{\bm{V}_s}^c$ is:
\begin{equation*}
  \left[L_{(d+1)m'\times (d+1)m'}^{s,(K)}\right]=\frac{1}{|K|}\left[b^{s,(K)}_{(d+1)m'\times (d+1)m'}\right],
\end{equation*}
where 
\begin{equation*}
    \left[b^{s,(K)}_{(d+1)m'\times ((d+1)m')}\right]_{i,(j-1)m' + r} = \int_{F_j}\psi_{r}^{(F_j)}\bvphi^{(K)}_{d(m-m')+i} \cdot \bn \, \mrmd\Gamma.
\end{equation*}


The local problem that depends on $\mu$ (equation \eqref{eqn:eq hdg method local uhat}) becomes the following matrix problem:
\begin{equation}\label{eqn:lin alg local prob mu stab-2}
  \begin{bmatrix}
    \left[I_{d(m-m')\times d(m-m')}\right] & \left[-D_{m \times d(m-m')}^{(K)}\right]^T\\
    \left[D^{(K)}_{m\times d(m-m')}\right] & \left[M^{(K),(L)}_{m \times m}\right]
  \end{bmatrix}
  \begin{bmatrix}
    \left[Q^{\mu,(K)}_{d(m-m')\times (d+1)m'}\right] \\
    \left[U^{\mu,(K)}_{m\times (d+1)m'}\right]
  \end{bmatrix}=
  \begin{bmatrix}
    -\left[b^{Q^{\mu},(K)}_{dm\times ((d+1)m')}\right]\\
    \left[b^{U^{\mu},(K)}_{m\times ((d+1)m')}\right]
  \end{bmatrix},
\end{equation}
where the mass-matrix from the stabilization term is:
\begin{equation*}
    \left[M^{s,(K)}_{m\times m}\right] = |K|\left[\left[L_{(d+1)m'\times (d+1)m'}^{s,(K)}\right]B^{(K)}_{(d+1)m'\times m}\right]^T \left[\left[L_{(d+1)m'\times (d+1)m'}^{s,(K)}\right]B^{(K)}_{(d+1)m'\times m}\right].
\end{equation*}
Then, similar to the Stab-1-HRT implementation, the above local matrix problem is solved using the
Cholesky factorization methodology after eliminating the matrix $\left[Q^{\mu,(K)}_{d(m-m')\times (d+1)m'}\right]$.

The local problem that depends on $f$ becomes:
{\small
\begin{equation}\label{eqn:lin alg local prob f stab-2}
  \begin{bmatrix}
    \left[I_{d(m-m')\times d(m-m')}\right] & \left[-D_{m \times d(m-m')}^{(K)}\right]^T\\
    \left[D^{(K)}_{m\times d(m-m')}\right] & \left[M^{(K),(L)}_{m \times m}\right]
  \end{bmatrix}
  \begin{bmatrix}
    \left[Q^{f,(K)}_{d(m-m')\times 1}\right] \\
    \left[U^{f,(K)}_{m\times 1}\right]
  \end{bmatrix}=
  \begin{bmatrix}
    -\left[0_{dm\times 1}\right]\\
    |K|\left[P^{(K)}f_{m\times 1}\right].
  \end{bmatrix}.
\end{equation}
}
The above matrix problem is also solved using the same Cholesky factors that were computed
while solving the above local problem that depended on $\mu$ alone.

\subsubsection{Global problem}

In this implementation, the element matrix and vector are computed as:
{\footnotesize
\begin{equation*}
\begin{split}
   \left[A^{(K)}_{((d+1)m') \times ((d+1)m')}\right] =& |K|\Bigg(\left[Q^{\mu,(K)}_{dm\times (d+1)m'}\right]^T \left[Q^{\mu,(K)}_{dm\times (d+1)m'}\right] + \\
   &\left[\left[L_{(d+1)m'\times (d+1)m'}^{s,(K)}\right]\left(B^{(K)}_{(d+1)m'\times m}\left[U^{\mu,(K)}_{m\times (d+1)m'}\right] - I_{(d+1)m' \times (d+1)m'}\right) \right]^T \\
   &\left[\left[L_{(d+1)m'\times (d+1)m'}^{s,(K)}\right]\left(B^{(K)}_{(d+1)m'\times m}\left[U^{\mu,(K)}_{m\times (d+1)m'}\right] - I_{(d+1)m' \times (d+1)m'}\right) \right]\Bigg), \\
   \left[b^{(K)}_{((d+1)m') \times 1}\right] =& \left[U^{\mu,(K)}_{m\times (d+1)m'}\right]^T \left[P^{(K)}f_{m\times 1}\right] |K|.
\end{split}
\end{equation*}
}
Then, similar to the Stab-1-HRT implementation, the above element matrices and element vectors are assembled to form the global problem
\begin{equation*}
  \left[A_{m'n_F \times m'n_{F}}\right] \left[\what{\mrm{u}}_{m'n_{F} \times 1}\right]=
  \left[ b_{m'n_{F} \times 1} \right]
\end{equation*}
and this global problem is solved using the sparse Cholesky factorization.

\section{Numerical results} \label{sec: numerical results}

We present results comparing the three implementations. 
We consider the Poisson problem with $f=8\pi^2 \sin(2\pi x_1)\sin(2\pi x_2)$
in the domain $(0,1)^2$. The domain is first split into 16 uniform quadrilateral
elements along each direction. Each quadrilateral element is further
split into two triangular elements. This leads to a uniform mesh of 512
triangular elements. Polynomial degrees 1 to 20 are considered. Zero-Dirichlet 
boundary condition is imposed on all the four sides of the domain. All there implementations
were first validated to make sure they yield identical solutions.

All numerical experiments were performed in MATLAB. We use a workstation with
Intel(R) Core(TM) i7-8700 processor. The processor has six cores and hyperthreading. 
We also note that MATLAB uses the multi-threaded MKL BLAS backend for certain matrix and 
vector manipulations. Hence, there is some inherent parallelism in our implementation. For 
the global problem solution, we use the sparse direct solver available in MATLAB. MATLAB 
uses the sparse multi-threaded Cholesky solver CHOLMOD \cite{chen2008algorithm} for our 
symmteric positive definite problem when performing
\texttt{x=A\textbackslash b}.

Tables \ref{tab: one time and local} and \ref{tab: global and total} show the time 
consumed by the one-time operations (that are performed just once in the reference element), 
the local problem solution in all the elements, global problem solution and 
the total solution time. The breakdown of the time taken by the different 
components of the local problem solution are shown in tables \ref{tab: additional RT basis and its divergence}
and \ref{tab: local matrix problem}. The percentage benefit of the new implementations 
compared to the Usual-HRT implementation is shown in table \ref{tab: percent benefit}.

\begin{table}[tb]
    \centering
{\tiny
\begin{tabular}{ccccccc}
    \hline
$k$&\multicolumn{3}{c|}{One-time operations} & \multicolumn{3}{c|}{Local problem solutions} \\ \hline
&	{\footnotesize Usual-HRT} &	{\footnotesize Stab-1-HRT} & {\footnotesize Stab-2-HRT} & {\footnotesize Usual-HRT} & {\footnotesize Stab-1-HRT} & {\footnotesize Stab-2-HRT} \\ \hline
1&8.54E-03&4.66E-03&6.14E-03&1.07E-01&8.29E-02&9.56E-02 \\ 
2&3.04E-03&2.11E-03&2.19E-03&1.65E-01&1.17E-01&1.26E-01 \\ 
3&3.46E-03&8.79E-04&8.62E-04&2.36E-01&1.78E-01&1.89E-01 \\ 
4&1.20E-03&1.20E-03&1.25E-03&3.59E-01&2.76E-01&2.90E-01 \\ 
5&1.77E-03&1.70E-03&1.83E-03&5.31E-01&4.23E-01&4.40E-01 \\ 
6&2.43E-03&2.52E-03&2.45E-03&8.82E-01&7.23E-01&6.70E-01 \\ 
7&3.37E-03&3.36E-03&3.39E-03&1.25E+00&1.05E+00&1.07E+00 \\ 
8&4.93E-03&5.01E-03&4.92E-03&1.74E+00&1.45E+00&1.47E+00 \\ 
9&6.85E-03&6.85E-03&7.16E-03&2.28E+00&1.99E+00&2.02E+00 \\ 
10&9.15E-03&9.17E-03&9.00E-03&3.07E+00&2.62E+00&2.68E+00 \\ 
11&1.21E-02&1.26E-02&1.25E-02&3.98E+00&3.43E+00&3.47E+00 \\ 
12&1.70E-02&1.72E-02&1.71E-02&5.13E+00&4.41E+00&4.59E+00 \\ 
13&2.33E-02&2.35E-02&2.37E-02&6.56E+00&5.83E+00&5.84E+00 \\ 
14&3.01E-02&2.99E-02&3.02E-02&8.27E+00&7.26E+00&7.39E+00 \\ 
15&4.01E-02&3.99E-02&4.04E-02&1.06E+01&9.26E+00&9.34E+00 \\ 
16&5.18E-02&5.12E-02&5.46E-02&1.30E+01&1.14E+01&1.16E+01 \\ 
17&7.16E-02&7.28E-02&7.09E-02&1.60E+01&1.40E+01&1.41E+01 \\ 
18&9.03E-02&8.97E-02&9.03E-02&1.95E+01&1.71E+01&1.73E+01 \\ 
19&1.25E-01&1.25E-01&1.26E-01&2.51E+01&2.07E+01&2.09E+01 \\ 
20&1.68E-01&1.63E-01&1.66E-01&2.99E+01&2.48E+01&2.51E+01 \\ \hline
    \end{tabular}}
    \caption{Comparison of the time (in seconds) for one-time operations and local problem solutions.}
    \label{tab: one time and local}
\end{table}

\begin{table}[tb]
    \centering
{\tiny
\begin{tabular}{ccccccc}
    \hline
$k$&\multicolumn{3}{c|}{Global problem solution} & \multicolumn{3}{c|}{Total solution} \\ \hline
&	{\footnotesize Usual-HRT} &	{\footnotesize Stab-1-HRT} & {\footnotesize Stab-2-HRT} & {\footnotesize Usual-HRT} & {\footnotesize Stab-1-HRT} & {\footnotesize Stab-2-HRT} \\ \hline
1&3.02E-02&3.00E-02&3.12E-02&1.46E-01&1.18E-01&1.33E-01 \\ 
2&2.76E-02&2.66E-02&3.06E-02&1.96E-01&1.46E-01&1.59E-01 \\ 
3&2.70E-02&2.65E-02&3.04E-02&2.66E-01&2.06E-01&2.20E-01 \\ 
4&3.03E-02&3.04E-02&3.41E-02&3.90E-01&3.08E-01&3.26E-01 \\ 
5&3.44E-02&3.35E-02&3.77E-02&5.67E-01&4.58E-01&4.79E-01 \\ 
6&3.96E-02&3.92E-02&4.19E-02&9.24E-01&7.65E-01&7.14E-01 \\ 
7&4.66E-02&4.65E-02&5.08E-02&1.30E+00&1.10E+00&1.13E+00 \\ 
8&5.24E-02&5.07E-02&5.54E-02&1.80E+00&1.51E+00&1.53E+00 \\ 
9&5.81E-02&5.84E-02&6.35E-02&2.35E+00&2.06E+00&2.09E+00 \\ 
10&6.54E-02&6.52E-02&6.99E-02&3.15E+00&2.69E+00&2.76E+00 \\ 
11&7.22E-02&7.10E-02&7.73E-02&4.06E+00&3.51E+00&3.56E+00 \\ 
12&8.00E-02&7.92E-02&8.51E-02&5.22E+00&4.51E+00&4.69E+00 \\ 
13&9.06E-02&8.95E-02&9.55E-02&6.68E+00&5.94E+00&5.96E+00 \\ 
14&1.00E-01&1.01E-01&1.05E-01&8.40E+00&7.39E+00&7.53E+00 \\ 
15&1.12E-01&1.10E-01&1.16E-01&1.07E+01&9.41E+00&9.50E+00 \\ 
16&1.25E-01&1.23E-01&1.29E-01&1.32E+01&1.16E+01&1.17E+01 \\ 
17&1.39E-01&1.37E-01&1.42E-01&1.62E+01&1.42E+01&1.43E+01 \\ 
18&1.52E-01&1.50E-01&1.55E-01&1.97E+01&1.73E+01&1.75E+01 \\ 
19&1.62E-01&1.60E-01&1.69E-01&2.54E+01&2.09E+01&2.12E+01 \\ 
20&1.81E-01&1.78E-01&1.85E-01&3.03E+01&2.51E+01&2.55E+01 \\ \hline
    \end{tabular}}
    \caption{Comparison of time (in seconds) for global problem solution and total solution.}
    \label{tab: global and total}
\end{table}

\begin{table}[t]
    \centering
    {\tiny
\begin{tabular}{ccccccc}
    \hline
$k$&\multicolumn{3}{c|}{Additional RT basis} & \multicolumn{3}{c|}{Div. matrix of additional RT basis} \\ \hline
&	{\footnotesize Usual-HRT} &	{\footnotesize Stab-1-HRT} & {\footnotesize Stab-2-HRT} & {\footnotesize Usual-HRT} & {\footnotesize Stab-1-HRT} & {\footnotesize Stab-2-HRT} \\ \hline
1&2.50E-02&2.52E-02&2.59E-02&2.74E-02&-&- \\ 
2&5.57E-02&5.35E-02&5.35E-02&4.43E-02&-&- \\ 
3&1.05E-01&1.05E-01&1.04E-01&5.99E-02&-&- \\ 
4&1.88E-01&1.88E-01&1.88E-01&8.57E-02&-&- \\ 
5&3.16E-01&3.18E-01&3.17E-01&1.15E-01&-&- \\ 
6&5.49E-01&5.48E-01&5.23E-01&1.61E-01&-&- \\ 
7&8.37E-01&8.37E-01&8.36E-01&2.08E-01&-&- \\ 
8&1.23E+00&1.21E+00&1.21E+00&2.72E-01&-&- \\ 
9&1.68E+00&1.71E+00&1.70E+00&3.35E-01&-&- \\ 
10&2.30E+00&2.29E+00&2.30E+00&4.29E-01&-&- \\ 
11&3.05E+00&3.05E+00&3.05E+00&5.31E-01&-&- \\ 
12&3.98E+00&3.97E+00&4.03E+00&6.83E-01&-&- \\ 
13&5.12E+00&5.24E+00&5.16E+00&8.89E-01&-&- \\ 
14&6.52E+00&6.57E+00&6.62E+00&1.08E+00&-&- \\ 
15&8.42E+00&8.47E+00&8.46E+00&1.37E+00&-&- \\ 
16&1.04E+01&1.05E+01&1.05E+01&1.61E+00&-&- \\ 
17&1.28E+01&1.29E+01&1.29E+01&2.05E+00&-&- \\ 
18&1.58E+01&1.58E+01&1.58E+01&2.41E+00&-&- \\ 
19&1.93E+01&1.92E+01&1.92E+01&4.31E+00&-&- \\ 
20&2.32E+01&2.31E+01&2.32E+01&5.04E+00&-&- \\ \hline
    \end{tabular}}
    \caption{Comparison of time (in seconds) for computation of the additional RT basis functions and their contribution to the divergence matrix for all the elements in the mesh.}
    \label{tab: additional RT basis and its divergence}
\end{table}

The new implementations Stab-1-HRT and Stab-2-HRT are faster than
the Usual-HRT implementation for all the polynomial degrees. From Table \ref{tab: percent benefit},
we observe that they are 10-20\% faster depending on the polynomial degree. 
Stab-1-HRT is slightly faster than Stab-2-HRT for nearly all polynomial 
degrees (except degree six).

The local problem solutions consume the majority of the total solution time. These local 
problem solution solutions are faster for the 
new Stab-1-HRT and Stab-2-HRT implementations compared to the usual-HRT implementations. 
From Tables 
\ref{tab: additional RT basis and its divergence}
and \ref{tab: local matrix problem} (which show the breakdown of the 
local problem solutions), we observe that this performance benefit is essentially
because in the new Stab-1-HRT and Stab-2-HRT implementations, we need not compute the 
derivative of the additional RT basis functions, i.e., $\bvphi_{dm+1},$$\dots,$$\bvphi_n$. 
and their contribution to the divergence matrix of the local problems. 
However, in the usual-HRT implementation, these derivatives and their contribution 
to the divergence matrix must be computed. In the Stab-1-HRT and Stab-2-HRT implementations, 
there is an overhead associated with the computation of the stabilization mappings. However, the benefits
from not computing the derivative of the additional RT basis functions and their contribution
to the divergence matrix significantly outweights this overhead. Hence, the new 
implementations of the hybridized Raviart-Thomas method yield significant (10-20\%) performance 
benefit compared to the usual implementation.
\begin{table}[b]
    \centering
    {\tiny
\begin{tabular}{ccccccc}
    \hline
$k$&\multicolumn{3}{c|}{Local matrix problem} \\ \hline
&	{\footnotesize Usual-HRT} &	{\footnotesize Stab-1-HRT} & {\footnotesize Stab-2-HRT} \\ \hline
1&3.47E-02&3.88E-02&4.98E-02 \\ 
2&4.04E-02&4.20E-02&5.07E-02 \\ 
3&4.27E-02&4.78E-02&5.86E-02 \\ 
4&4.93E-02&5.42E-02&6.88E-02 \\ 
5&5.69E-02&6.34E-02&8.10E-02 \\ 
6&1.15E-01&1.21E-01&9.54E-02 \\ 
7&1.37E-01&1.45E-01&1.68E-01 \\ 
8&1.58E-01&1.58E-01&1.85E-01 \\ 
9&1.65E-01&1.82E-01&2.21E-01 \\ 
10&2.06E-01&2.03E-01&2.54E-01 \\ 
11&2.36E-01&2.34E-01&2.82E-01 \\ 
12&2.74E-01&2.69E-01&3.40E-01 \\ 
13&3.27E-01&3.70E-01&4.62E-01 \\ 
14&3.95E-01&4.26E-01&5.03E-01 \\ 
15&4.55E-01&4.87E-01&5.74E-01 \\ 
16&5.54E-01&5.44E-01&6.88E-01 \\ 
17&6.07E-01&6.41E-01&7.81E-01 \\ 
18&7.19E-01&7.51E-01&9.29E-01 \\ 
19&8.61E-01&9.04E-01&1.12E+00 \\ 
20&1.02E+00&1.05E+00&1.31E+00 \\ \hline
    \end{tabular}}
    \caption{Comparison of time (in seconds) for the solution of the local matrix problems using the Cholesky decomposition for each element in the mesh.}
    \label{tab: local matrix problem}
\end{table}

\begin{table}[bt]
    \centering
    {\tiny
\begin{tabular}{ccccccc}
    \hline
$k$&\multicolumn{2}{c|}{Local problem solution} & \multicolumn{2}{c|}{Total solution} \\ \hline
&	{\footnotesize Stab-1-HRT} & {\footnotesize Stab-2-HRT} &  {\footnotesize Stab-1-HRT} & {\footnotesize Stab-2-HRT} \\ \hline
1&22.65&10.90&19.40&8.93 \\ 
2&28.84&23.40&25.32&18.63 \\ 
3&24.35&20.03&22.72&17.46 \\ 
4&23.04&19.06&21.16&16.55 \\ 
5&20.34&17.23&19.21&15.53 \\ 
6&17.94&23.99&17.15&22.64 \\ 
7&16.20&14.47&15.58&13.59 \\ 
8&16.85&15.67&16.41&15.01 \\ 
9&12.77&11.44&12.40&10.88 \\ 
10&14.87&12.83&14.53&12.38 \\ 
11&13.85&12.71&13.58&12.32 \\ 
12&13.95&10.49&13.70&10.20 \\ 
13&11.15&11.01&10.97&10.74 \\ 
14&12.21&10.64&12.01&10.42 \\ 
15&12.31&11.56&12.16&11.36 \\ 
16&12.50&11.24&12.36&11.04 \\ 
17&12.40&11.64&12.24&11.47 \\ 
18&12.20&11.08&12.07&10.94 \\ 
19&17.59&16.47&17.40&16.25 \\ 
20&17.15&16.01&16.98&15.82 \\ \hline
    \end{tabular}}
    \caption{Percentage performance benefit of the Stab-1-HRT and Stab-2-HRT implementations over the Usual-HRT implementation.}
    \label{tab: percent benefit}
\end{table}

\clearpage

\section{Concluding remarks and ongoing work}
As pointed out in \cite{bcockburn2023hybridizable}, although the choice of the space $\bm{V}_a(K)$ is not unique, as we have also seen here, the {\em smallest} of these spaces, $\nabla W(K)=\nabla {P}_k(K)$, is actually unique. It remains to be seen it we still retain an advantage over the usual  implementation of the hybridized Raviart-Thomas method for this choice. This constitutes the subject of ongoing work.

The choice of the Raviart-Thomas method as the mixed method is by no means restrictive, as the approach displayed here can also be applied to any other mixed method, standard or generated by the theory of M-decompositions; see the review \cite{CockburnFuShi18} and the references therein. Also,
the application of this approach to general HDG methods
can be carried out very easily, as we are going to show elsewhere.

Finally, note that, although we have exploited the fact that the tensor-valued function $c$ is the identity for our  model second-order elliptic problem, it is easy to extend what has been done to a general elliptic problem. This is also part of our ongoing work.

\bibliographystyle{plain}
\bibliography{papers}

\end{document}